\documentstyle[amsmath,thm,pb-diagram,amssymb,leqno,epic,graphicx,xy,undertilde,amssymb,epigraph]{article}
\xyoption{matrix}
\xyoption{arrow}
\newfont\got{eufm10}

\newtheorem{proposition}{Proposition}[section]
\newtheorem{thm}[proposition]{Theorem}
\newtheorem{cor}[proposition]{Corollary}
\newtheorem{lemma}[proposition]{Lemma}

\newtheorem{claim}[proposition]{Claim}

\newcounter{secnum}

\begin{document}

\begin{center}
{\Large \bf Definable Hamel bases and ${\sf AC}_\omega({\mathbb R})$}
\end{center}

\begin{center}
{\bf Vladimir Kanovei\footnote{The first author gratefully acknowledges support through the RFBR grant 17-01-00705.}}\\
IITP RAS\\
Bolshoy Karetny per. 19, build.1\\
Moscow 127051 Russia
\end{center}
\begin{center}
{\bf Ralf Schindler\footnote{Funded by the Deutsche Forschungsgemeinschaft (DFG, German Research Foundation) under Germany's Excellence Strategy
EXC 2044 –390685587, Mathematics Münster: Dynamics - Geometry - Structure.}}\\
Institut f\"ur Mathematische Logik und Grundlagenforschung,
Universit\"at M\"unster\\
Einsteinstr.~62,
48149 M\"unster, Germany\\
\end{center}


\begin{abstract} There is a model of ${\sf ZF}$ with a $\Delta^1_3$ definable
Hamel bases in which ${\sf AC}_\omega({\mathbb R})$ fails.
\end{abstract}

Answering a question from \cite[p.\ 433]{pincus_prikry} it was shown in \cite{BSchWY} that there is a Hamel basis in the Cohen--Halpern--L\'evy model. In this paper we show that in a variant of this model, there is a projective, in fact $\Delta^1_3$, Hamel basis.

Throughout this paper, by a {\em Hamel basis} we always mean a basis for ${\mathbb R}$, construed as
a vector space over ${\mathbb Q}$. We denote by $E_0$ the {\em Vitali equivalence
relation}, $x E_0 y$ iff $x-y \in {\mathbb Q}$ for $x$, $y \in {\mathbb R}$.
We also write $[x]_{E_0} = \{ y \colon y E_0 x \}$ for the $E_0$--equivalence class
of $x$.
A transversal for the set of all $E_0$--equivalence classes picks exactly
one member from each $[x]_{E_0}$. The range of any such transversal is also called a
{\em Vitali set}. If we identify ${\mathbb R}$ with the Cantor space ${}^\omega 2$,
then $x E_0 y$ iff $\{ n \colon x(n) \not= y(n) \}$ is finite.

A set $\Lambda \subset {\mathbb R}$ is a {\em Luzin set} iff $\Lambda$ is uncountable but
$\Lambda \cap M$ is at most countable for every meager set $M \subset {\mathbb R}$.
A set $S \subset {\mathbb R}$ is a {\em Sierpi\'nski set} iff $S$ is uncountable but
$S \cap N$ is at most countable for every null set $N \subset {\mathbb R}$
(``null'' in the sense of Lebesgue measure). A set $B \subset {\mathbb R}$
is a {\em Bernstein set} iff $B \cap P \not= \emptyset \not= P \setminus B$ for every
perfect set $P \subset {\mathbb R}$. A {\em Burstin basis} is a Hamel basis which is also a Bernstein set. It is easy to see that $B \subset {\mathbb R}$ is a Burstin basis 
iff $B$ is a Hamel basis and $B \cap P \not= \emptyset$ for every perfect $P
\subset {\mathbb R}$.

A set $m \subset {\mathbb R} \times {\mathbb R}$ is called a {\em Mazurkiewicz
set} iff $m \cap \ell$ for every straight line $\ell \subset {\mathbb R} \times {\mathbb R}$.

By ${\sf AC}_\omega({\mathbb R})$ we mean the statement that for all
sequences $(A_n \colon n<\omega)$ such that $\emptyset \not= A_n \subset {\mathbb R}$ for all
$n<\omega$ there is some choice function $f \colon \omega \rightarrow {\mathbb R}$, i.e.,
$f(n) \in A_n$ for all $n<\omega$.

D.\ Pincus and K.\ Prikry 
study the Cohen-Halpern-L\'evy model $H$ in \cite{pincus_prikry}. The model
$H$ 
is obtained by adding a countable set of Cohen reals (say over $L$) without adding their enumeration;
$H$
does not satisfy ${\sf AC}_\omega({\mathbb R})$. It is shown in \cite{pincus_prikry} that there is
a Luzin set in $H$, so that in ${\sf ZF}$, the existence of a 
Luzin set does not even imply ${\sf AC}_\omega({\mathbb R})$.
\cite[Theorems 1.7 and 2.1]{BSchWY} show that in $H$ there is a Bernstein set as
well as a Hamel basis. As in ${\sf ZF}$ the existence of a Hamel basis implies the existence
of a Vitali set, the latter also reproves Feferman's result (see \cite{pincus_prikry})
according to which there is a Vitali set in $H$.


Therefore, in ${\sf ZF}$ the conjunction of the following statements (1), (3), and (5) 
(which in ${\sf ZF}$ implies (4)) does
not yield ${\sf AC}_\omega({\mathbb R})$.   

\begin{itemize}\itemsep0.0pt
\item[(1)] There is a Luzin set.
\item[(2)] There is a Sierpi\'nski set. 
\item[(3)] There is a Bernstein set.
\item[(4)] There is a Vitali set.
\item[(5)] There is a Hamel basis.
\item[(6)] There is a Burstin basis.
\item[(7)] There is a Mazurciewicz set.
\end{itemize}

(2) is false in $H$, see \cite[Lemma 1.6]{BSchWY}. We neither know if (6) is true
in $H$, nor do we know if (7) is true in $H$.
We aim to prove that in ${\sf ZF}$, the conjunction of {\em all} of these statements
does not imply ${\sf AC}_\omega({\mathbb R})$, even if the respective sets
are required to be projective.

The Luzin set which \cite[Theorem on p.\ 429]{pincus_prikry} constructs is $\Delta^1_2$. 
In ${\sf ZFC}$, there is no analytic Hamel basis (see \cite{sierpinski},\cite{sierpinski_2}, \cite{jones}), but by a theorem of
A.\ Miller, in $L$ there is a coanalytic Hamel basis, see \cite[Theorem 9.26]{a_miller}; see also e.g.\ \cite[Corollary 2 and Lemma 4]{talks-kiel}. On the other hand, it can be verified that the model from \cite{BSchWY}
doesn't have a projective Vitali set.\footnote{To display our ignorance: we don't know if the model from \cite{brendle_et_al} has a definable Hamel basis.}
For the convenience of the reader as well as to motivate what is to come, we shall
sketch the proof of this at the beginning of the first section,
see Lemma \ref{motivation}.

The papers \cite{brendle_et_al} and \cite{mazurkiewicz} produce models of
{\sf ZF} plus {\sf DC} plus (6) and {\sf ZF} plus {\sf DC} plus (7), respectively.
By another theorem of A.\ Miller, see \cite[Theorem 7.21]{a_miller}, in $L$ there is
a coanalytic Mazurkiewicz set. It is not known if there is a Mazurkiewicz set which is Borel.

The result of the current paper is the following.

\begin{thm}\label{main-thm} There is a model of ${\sf ZF}$ plus $\lnot {\sf AC}_\omega({\mathbb R})$ in which the following hold true.
\begin{itemize}\itemsep0.0pt
\item[(a)] There is a $\Delta^1_2$ Luzin set.
\item[(b)] There is a $\Delta^1_2$ Sierpi\'nski set. 
\item[(c)] There is a $\Delta^1_3$ Bernstein set.
\item[(d)] There is a $\Delta^1_3$ Hamel basis.
\end{itemize}
\end{thm} 

\section{Jensen's perfect set forcing, revisited.}

In what follows, we shall mostly think of reals as elements of the Cantor
space ${}^\omega 2$.
We shall need a variant of the Cohen-Halpern-L\'evy model.

In order to get a definable Hamel basis in the absence of ${\sf AC}_\omega({\mathbb R})$ forces us to indeed work
with a model which is different from the original 
Cohen-Halpern-L\'evy model. This follows from the followig folklore result which we include
here as a motvation for what is to come. Recall, see \cite[Lemma 1.1]{BSchWY},
that a Hamel basis trivially produces a Vitali set.

Recall that the original 
Cohen-Halpern-L\'evy model is produced as follows. See \cite{pincus_prikry}, see also \cite[p.\ 3567]{BSchWY}. Let $g$ be ${\mathbb C}(\omega)$-generic over $L$, and let $A$ denote the countable set of Cohen reals which
${\mathbb C}(\omega)$ adds. Then 
\begin{eqnarray}\label{original}
H = {\sf HOD}^{L[g]}_{A \cup \{ A \} }.
\end{eqnarray}

\begin{lemma}\label{motivation} The Cohen-Halpern-L\'evy model from (\ref{original}) does not have a Vitali set which is definable in $H$ from ordinals and reals.\footnote{Recall, though, that
$H$ does have a Hamel basis which is definable from the {\em set} $A$ of Cohen reals, see the proof in \cite[section 2]{BSchWY}.}
\end{lemma}

{\em Proof.} Let $g$ be ${\mathbb C}(\omega)$-generic over $L$. It suffices to prove that there is no $a \in [A]^{<\omega}$ such that in $L[g]$ there is a Vitali set in ${}^\omega 2$ which is
definable from ordinals and $a$. 

Suppose otherwise. By minimizing the ordinal parameters, we may fix
$a \in [A]^{<\omega}$ such that in $L[g]$ there is a Vitali set which is
definable $a$, say via the formula $\varphi(-,a)$. Let $c \in A \setminus a$, and say
$n < \omega$ and $s \in {}^n 2$ are such that 
\begin{eqnarray}\label{f1}
L[g] \models \varphi(s^\frown c \upharpoonright [n,\omega),a).
\end{eqnarray}
Let ${\dot c}$ be a canonical ${\mathbb C}(\omega)$-name for $c$, 
in particular, ${\dot c}^g = c$, let ${\dot a}$ be a canonical ${\mathbb C}(\omega)$-name for $a$, in particular, ${\dot a}^g = a$,
and pick $p \in g$ such that 
\begin{eqnarray}\label{f2}
p \Vdash_L^{{\mathbb C}(\omega)}
\varphi({\check s}^\frown {\dot c} \upharpoonright [{\check n},\omega),{\dot a}).
\end{eqnarray}
Let $g^*$ be ${\mathbb C}(\omega)$-generic over $L$ which is identical with $g$
except for that $g^*$ incorporates a finite nontrivial variant of $g$ only in the coordinate of ${\mathbb C}(\omega)$ which
gives rise to $c$ in such a way that ${\dot c}^{g^*} \upharpoonright [n,\omega)
\not= c \upharpoonright [n,\omega)$, but ${\dot c}^{g^*}$ is $E_0$-equivalent with $c$.
We have that $L[g^*]=L[g]$ and ${\dot a}^{g^*} = a$, and (\ref{f2}) yields that
\begin{eqnarray}\label{f3}
L[g] \models \varphi(s^\frown {\dot c}^{g^*} \upharpoonright [n,\omega),a).
\end{eqnarray}
(\ref{f1}) and (\ref{f3}) contradictthe fact that $\varphi(-,a)$ defines a Vitali set
in $L[g]$. $\square$

\bigskip
The same argument shows that the model from \cite{brendle_et_al} doesn't have a Vitali
set which is definable from ordinal and real parameters.

In order to construct our model, we now need to introduce a variant of 
Jensen's variant of Sacks forcing, see \cite{jensen} (see also \cite[Definition 6.1]{kanovei-lyubetsky}), which
we shall call ${\mathbb P}$. The reason why we can't work with Jensen's forcing
directly is that it does not seem to have the Sacks property (see e.g.\ \cite[Definition 2.15]{brendle_et_al}).

By way of notation, if ${\mathbb Q}$ is a forcing and $N>0$ is any ordinal, then
${\mathbb Q}(N)$ denotes the finite support product of $N$ copies of ${\mathbb Q}$,
ordered component-wise.
In this paper, we shall only consider ${\mathbb Q}(N)$ for $N \leq \omega$.
If $\alpha$ is a limit ordinal, then $<_{J_\alpha}$ denotes the canonical
well-ordering of $J_\alpha$, see \cite[Definition 5.14 and p.\ 79]{book},\footnote{The reader unfamiliar with the $J$-hierarchy may read $L_\alpha$ instead of $J_\alpha$.} and $<_L = \bigcup \,
\{ <_{J_\alpha} \colon \alpha$ is a limit ordinal $\}$.

Let us work in $L$ until further notice. Let us first define $(\alpha_\xi , \beta_\xi
\colon \xi < \omega_1)$ as follows: $\alpha_\xi =$ the least $\alpha > {\rm sup}(\{
\beta_{\bar \xi} \colon {\bar \xi} < \xi \})$ such that $J_\alpha \models 
{\sf ZFC}^-$,\footnote{Here, ${\sf ZFC}^-$ denotes ${\sf ZFC}$ without the power
set axiom. Every $J_\alpha$ satisfies the strong form of ${\sf AC}$ according to which
every set is the surjective image of some 
ordinal. In the absence of $V=L$, one has to be careful about how to formulate ${\sf ZFC}^-$, see \cite{gitman}.} and $\beta_\xi =$ the least $\beta > \alpha_\xi$ such that $\rho_\omega(J_\beta)=\omega$ (see \cite[Definition 11.22]{book}; $\rho_\omega(J_\beta)=\omega$ is
equivalent with ${\cal P}(\omega) \cap J_{\beta+\omega} \not\subset J_\beta$). 

We shall also make use of a sequence 
$(f_\xi \colon \xi<\omega)$ which is defined as follows. 
Let $({\bar f}_\xi \colon \xi<\omega)$ be defined by the following trivial recursion:
${\bar f}_\xi$ be the $<_L$-least $f$
such that $f \in ({}^\omega J_{\omega_1} \cap J_{\omega_1}) 
\setminus \{ {\bar f}_{\bar \xi} \colon {\bar \xi}<\xi \}$. 
Then if $\pi$ denotes the G\"odel pairing function, see \cite[p.\ 35]{book},
we let $f_{\pi((\xi_1,\xi_2))}={\bar f}_{\xi_1}$.
We will then have that $f_\xi \in J_{\alpha_\xi}$ for all $\xi$, and
for each $f \in ({}^\omega J_{\omega_1} \cap J_{\omega_1})$
the set of $\xi$ such that $f=f_\xi$ is cofinal in $\omega_1$.

Let us then define $({\mathbb P}_\xi , {\mathbb Q}_\xi \colon \xi \leq \omega_1)$.
Each ${\mathbb P}_\xi$ will consist of perfect trees $T \subset {}^{<\omega} 2$ such
that if $T \in {\mathbb P}_\xi$ and $s \in T$, then $T_s = \{ t \in T \colon
t \subset s \vee s \subset t \} \in {\mathbb P}_\xi$ as well.\footnote{We denote by $x \subset y$ the fact that $x$ is a (not necessarily proper) subset of $y$.}
Each ${\mathbb P}_\xi$ will be construed as a p.o.\ by stipulating $T \leq T'$ ($T$ ``is stronger than'' $T'$) iff
$T \subset T'$. We will have that ${\mathbb P}_\xi \in J_{\alpha_\xi}$ and
${\mathbb P}_{\bar \xi} \subset {\mathbb P}_\xi$
whenever ${\bar \xi} \leq \xi \leq \omega_1$.

To start with, let ${\mathbb P}_0$ be the set of all basic clopen sets $U_s = \{ t \in
{}^{<\omega} 2 \colon t \subset s \vee s \subset t \}$, where $s \in {}^\omega 2$.
If $\lambda \leq \omega_1$ is a limit ordinal, then ${\mathbb P}_\lambda =
\bigcup \{ {\mathbb P}_\xi \colon \xi<\lambda \}$. 

Now fix $\xi<\omega_1$, and suppose that ${\mathbb P}_\xi$ has already been
defined. We shall define ${\mathbb Q}_\xi$ and ${\mathbb P}_{\xi+1}$.

Let $g_\xi \in {}^\omega J_{\alpha_\xi}$ be the following $\omega$-sequence.
If there is some $N<\omega$ such that $f_\xi$ is an $\omega$-sequence of subsets of ${\mathbb P}_{\xi}(N)$,
each of which is predense in ${\mathbb P}_{\xi}(N)$,
then for each $n<\omega$ let $g_\xi(n)$ be the open dense set $$\{ (T_1, \ldots , T_N) \in {\mathbb P}_\xi(N) \colon \exists (T'_1 , \ldots , T'_N) \in f_\xi(n) \, (T_1, \ldots , T_N) \leq (T'_1 , \ldots , T'_N) \}{\rm , }$$ and write $N_\xi=N$.
Otherwise we just set $g_\xi(n) = {\mathbb P}_\xi(1)$ for each $n<\omega$, and write
$N_\xi=1$.
Let $d_\xi$ be the $<_{J_{\beta_\xi+\omega}}$-least $d \in {}^{\omega \times
\omega}({\cal P}({\mathbb 
P}_\xi) \cap J_{\alpha_\xi}) \cap J_{\beta_\xi+\omega}$ such that
\begin{enumerate}
\item[(i)] for each $(n,N) \in \omega \times \omega$, $d(n,N)$ is an open dense subset of ${\mathbb P}_\xi(N)$ which exists in $J_{\beta_\xi}$,
\item[(ii)] for each $N<\omega$ and each open dense subset $D$ of ${\mathbb P}_\xi(N)$ which exists in $J_{\beta_\xi}$ there is some $n<\omega$ with $d(n,N) \subset D$, 
\item[(iii)] $d(n,N_\xi) \subset g_\xi(n)$ for each $n<\omega$, and
\item[(iv)] $d(n+1,N) \subset d(n,N)$ for each $(n,N) \in \omega \times \omega$.
\end{enumerate}

Let us now look at the collection of all systems $(T^m_s \colon m<\omega , s \in {}^{<\omega} 2 )$ with the following properties.
\begin{enumerate}
\item[(a)] $T^m_s \in {\mathbb P}_\xi$ for all $m$, $s$,
\item[(b)] for each $T \in {\mathbb P}_\xi$ there are infinitely many $m<\omega$ with
$T^m_\emptyset = T$,
\item[(c)] $T^m_t \leq T^m_s$ for all $m$, $t \supset s$,
\item[(d)] ${\rm stem}(T^m_{s^\frown 0})$ and ${\rm stem}(T^m_{s^\frown 1})$
are incompatible elements of $T^m_s$ for all $m$, $s$, 
\item[(e)] if $(m,s) \neq (m',s')$, where $m$,$m' <n$ and ${\rm lh}(s)={\rm lh}(s') =n+1$ for some $n$, then
${\rm stem}(T^m_s)$ and ${\rm stem}(T^{m'}_{s'})$ are incompatible, and
\item[(f)] for all $N \leq n < \omega$ and all pairwise different $(m_1,s_1)$, $\ldots$, $(m_N,s_N)$ with $m_1$, $\ldots$, $m_N
< n$ and $s_1$, $\ldots$, $s_N \in {}^{n+1} 2$, $$(T^{m_1}_{s_1} , \ldots
, T^{m_N}_{s_N}) \in d_\xi(n,N).$$
\end{enumerate}
It is easy to work in $J_{\beta_\xi+\omega}$ and construct initial segments 
$(T^m_s \colon m<\omega , s \in {}^{<\omega} 2 , {\rm lh}(s) \leq n)$ of such a
system by induction on $n<\omega$. Notice that (f) formulates a constraint 
only for $m_1$, $\ldots$, $m_N
< {\rm lh}(s_1)-1 = \ldots = {\rm lh}(s_N)-1$, and writing $n={\rm lh}(s_1)-1$, there are $\sum_{N=1}^{n} \frac{(n \cdot 2^{n+1})!}{(n \cdot 2^{n+1} -N)!}$ (i.e.,
finitely many) such constraints. 

We let $(T^m_{s,\xi} \colon m<\omega , s \in {}^{<\omega} 2 )$ be the
$<_{\beta_\xi+\omega}$-least such system $(T^m_s \colon m<\omega , s \in {}^{<\omega} 2 )$.
For every $m<\omega$, $s \in {}^{<\omega} 2$, we let $$A^m_{s,\xi} = \bigcap_{n \geq {\rm lh}(s)} ( \bigcup_{\substack{t \supset s \\ {\rm lh}(t)=n}} \, T^m_t )
= \{ {\rm stem}(T_{t,\xi}^m)\upharpoonright k \colon t \supset s , k<\omega \}.$$
Notice that (e) implies that 
\begin{eqnarray}\label{incompatible}
A^m_{s,\xi} \cap A^{m'}_{s',\xi} \mbox{ is finite, unless } m=m' \mbox{ and }
s \subset s' \mbox{ or } s' \subset s.
\end{eqnarray}
(\ref{incompatible}) will imply that $A^m_{s,\xi}$ and $A^{m'}_{s',\xi}$ will be
incompatible in every ${\mathbb P}_\eta$, $\eta>\xi$, unless $m=m'$ and $s \subset s'$ or $s' \subset s$.

We set ${\mathbb Q}_\xi = \{ A^m_{s,\xi} \colon m<\omega , s \in {}^{<\omega} 2 \}$. Finally, we set ${\mathbb P}_{\xi+1} = {\mathbb P}_\xi \cup {\mathbb Q}_\xi$.

\begin{lemma}\label{1} Let $N<\omega$, $\xi < \omega_1$. $$D=\{ (T_1 , \ldots , T_N) \in {\mathbb Q}_\xi(N) \colon \mbox{stem}(T_i) \perp \mbox{stem}(T_j) \mbox{ for } i \neq j \}$$ is dense in ${\mathbb P}_{\xi+1}(N)$.\footnote{Here,
$\mbox{stem}(T_i) \perp \mbox{stem}(T_j)$ means that the stem of $T_i$ is incompatible
with the stem of $T_j$.}
\end{lemma}
{\em Proof.} Let $(T_1, \ldots , T_N) \in {\mathbb P}_{\xi+1}(N)$. For $i \in
\{ 1 , \ldots , N \}$ such that $T_i \in {\mathbb P}_\xi$ pick some $m_i<\omega$
such that $T_i = T^{m_i}_{\emptyset,\xi}$, and write $s_i = \emptyset$. This is possible by (b). If $i \in
\{ 1 , \ldots , N \}$ is such that $T_i \in {\mathbb Q}_\xi$, then say $T_i =
A^{m_i}_{s_i,\xi}$. Now pick $n>{\rm max}(\{ 
m_1 , \ldots , m_N \})$ and $t_1 \supset s_1$, $\ldots$, $t_N \supset s_N$ such that ${\rm lh}(t_1)
= \ldots = {\rm lh}(t_N) = n+1$ and the $(m_i,t_i)$ are pairwise different.

Then by (e) the finite sequences ${\rm stem}(T^{m_i}_{t_i,\xi})$ are pairwise incompatible, so that by $A^{m_i}_{t_i,\xi} \leq T^{m_i}_{t_i,\xi}$, the $A^{m_i}_{t_i,\xi}$ are pairwise incompatible. But
then $(A^{m_1}_{t_i,\xi}, \ldots , A^{m_N}_{t_N,\xi}) \in D$ and 
$(A^{m_1}_{t_i,\xi}, \ldots , A^{m_N}_{t_N,\xi}) \leq (T_1, \ldots , T_N)$. $\square$ 

\begin{lemma}\label{2} {\bf (Sealing)} Let $N<\omega$, $\xi < \omega_1$. If $D \in J_{\beta_\xi}$ is predense in ${\mathbb P}_\xi(N)$, then $D$ is predense in all ${\mathbb P}_{\eta}(N)$, $\eta \geq \xi$, $\eta \leq \omega_1$.
\end{lemma}
{\em Proof} by induction on $\eta$. The cases $\eta=\xi$ and $\eta$ being a limit
ordinal are trivial. Suppose $\eta \geq \xi$, $\eta < \omega_1$, and $D$ is predense in ${\mathbb P}_\eta(N)$. Write $D' = \{ (T_1, \ldots , T_N) \in {\mathbb P}_\eta(N) \colon \exists (T'_1 , \ldots , T'_N) \in D \, (T_1, \ldots , T_N) \leq (T'_1 , \ldots , T'_N) \}$.
As $\beta_\xi \leq \beta_\eta$, $D' \in J_{\beta_\xi}$ and by (ii) and (iv) there is some
$n_0<\omega$ with $d_\eta(n,N) \subset D'$ for every $n > n_0$.

To show that $D'$ (and hence $D$) is predense in ${\mathbb P}_{\eta+1}(N)$, by Lemma \ref{1}
it suffices to show that for all $(T_1, \ldots , T_N) \in {\mathbb Q}_\eta(N)$
there is some $(T'_1 , \ldots , T'_N) \in {\mathbb Q}_\eta(N)$,  $(T'_1 , \ldots , T'_N)
\leq (T_1, \ldots , T_N)$, and $(T'_1 , \ldots , T'_N)$ is below some element of $D'$.

So let $(A^{m_1}_{s_1,\eta}, \ldots , A^{m_N}_{s_N,\eta}) \in {\mathbb Q}_\eta(N)$
be arbitrary. Let $$n > {\rm max}( \{ n_0 , N-1 , m_1 , \ldots , m_N ,
{\rm lh}(s_1) , \ldots , {\rm lh}(s_N) \} ){\rm , }$$ and let $t_1 \supset s_1$, $\ldots$, $t_N \supset s_N$ be such that ${\rm lh}(t_1) = \ldots = {\rm lh}(t_N) = n+1$.
By increasing $n$ further if necessary, we may certainly assume that $t_1$, $\ldots$, $t_N$ are picked in such a way that $(m_1,t_1)$, $\ldots$, $(m_N,t_N)$ are pairwise different.
Then 
$$(T^{m_1}_{t_1,\eta}, \ldots , T^{m_N}_{t_N,\eta}) \in d_\eta(n,N) \subset D'$$
by (f). But $$(A^{m_1}_{t_1,\eta}, \ldots , A^{m_N}_{t_N,\eta}) \leq (T^{m_1}_{t_1,\eta}, \ldots , T^{m_N}_{t_N,\eta}){\rm , }$$ and also 
$$(A^{m_1}_{t_1,\eta}, \ldots , A^{m_N}_{t_N,\eta}) \leq
(A^{m_1}_{s_1,\eta}, \ldots , A^{m_N}_{s_N,\eta}){\rm , }$$ 
which means that $(A^{m_1}_{s_1,\eta}, \ldots , A^{m_N}_{s_N,\eta})$ is
compatible with
an element of $D'$. $\square$

\begin{cor}\label{3} Let $N<\omega$, $\xi < \omega_1$. $$\{ (T_1 , \ldots , T_N) \in {\mathbb Q}_\xi(N) \colon \mbox{stem}(T_i) \perp \mbox{stem}(T_j) \mbox{ for } i \neq j \}$$ is predense in ${\mathbb P}(N)$.
\end{cor}

\begin{lemma}\label{4}
Let $N<\omega$. ${\mathbb P}(N)$ has the c.c.c.
\end{lemma}
{\em Proof.} Let $A \subset {\mathbb P}(N)$ be a maximal antichain, $A \in L$.
Let $j \colon J_\beta \rightarrow J_{\omega_2}$ be elementary and such that $\beta
< \omega_1$ and $\{ {\mathbb P} , A \} \subset {\rm ran}(j)$. Write $\xi = {\rm crit}(j)$. We have that
$j^{-1}({\mathbb P}(N))={\mathbb P}(N) \cap J_\xi={\mathbb P}_\xi(N)$ and $j^{-1}(A)=A \cap J_\xi = A \cap {\mathbb P}_\xi(N) \in J_\beta$ is a maximal antichain in ${\mathbb P}_\xi(N)$. Moreover, $\beta_\xi > \beta$, so that by Lemma \ref{3} $A \cap {\mathbb P}_\xi(N)$ is predense in ${\mathbb P}(N)$. This means that $A = A \cap
{\mathbb P}_\xi$ is countable. $\square$

\begin{lemma}\label{5} Let $N<\omega$. $(c_1 , \ldots c_N) \in {}^N ({}^\omega 2)$ is
${\mathbb P}(N)$-generic over $L$ iff for all $\xi<\omega_1$ there 
is an injection $t \colon \{ 1 , \ldots , N \} \rightarrow {\mathbb Q}_\xi$ such that
for all $i \in \{ 1 , \ldots , N \}$, $c_i \in [t(i)]$.
\end{lemma}
{\em Proof.} ``$\Longrightarrow$'': This readily follows from Corollary \ref{3}.

``$\Longleftarrow$'': Let $A \subset {\mathbb P}(N)$ be a maximal antichain, $A \in L$. By Lemma \ref{4}, we may certainly pick some $\xi<\omega_1$ with $A \subset {\mathbb 
P}_\xi(N)$ and $A \in J_{\alpha_\xi}$. Say $n_0$ is such that $d_\xi(n,N) \subset
\{ (T_1, \ldots , T_N) \in {\mathbb P}_\xi \colon \exists (T'_1, \ldots , T'_N) \in A \, (T_1, \ldots , T_N) \leq (T'_1, \ldots , T'_N) \}$ for all $n \geq n_0$.
By our hypothesis, we may pick pairwise different $(m_1,s_1)$, $\ldots$, $(m_N,s_N)$
with ${\rm lh}(s_1) = \ldots = {\rm lh}(s_N)=n+1$ for some $n \geq n_0$ and
$c_i \in [T^{m_i}_{s_i},\xi]$ for all $i \in \{ 1, \ldots , N \}$.
But then $(T^{m_1}_{s_i,\xi}, \ldots , T^{m_N}_{s_N})$ is below an element of $A$,
which means that the generic filter given by $(c_1, \ldots , c_N)$ meets $A$. $\square$

\begin{cor} Let $N<\omega$, and let $(c_1 , \ldots c_N) \in {}^N ({}^\omega 2)$ be
${\mathbb P}(N)$-generic over $L$. If $x \in L[(c_1 , \ldots c_N)]$ is ${\mathbb P}$-generic over $L$, then $x \in \{ c_1 , \ldots c_N \}$. 
\end{cor}

{\em Proof.} If $x \in L[(c_1 , \ldots c_N)]$ is ${\mathbb P}$-generic over $L$, then 
$(c_1 , \ldots c_N, x) \in {}^{N+1} ({}^\omega 2)$ is
${\mathbb P}(N+1)$-generic over $L$, hence $x \notin L[(c_1 , \ldots c_N)]$.
Contradiction! $\square$

\begin{cor}\label{definability}
Let $N<\omega$, and let $(c_1 , \ldots c_N) \in {}^N ({}^\omega 2)$ be
${\mathbb P}(N)$-generic over $L$. Then inside $L[(c_1 , \ldots c_N)]$, $\{ c_1 , \ldots c_N \}$ is a (lightface) $\Pi^1_2$ set.
\end{cor}

{\em Proof.} Let $\varphi(x)$ express that for all $\xi<\omega_1$ there 
is some $T \in {\mathbb Q}_\xi$ such that
$x \in [T]$. The formula $\varphi(x)$ may be written in a $\Pi^1_2$ fashion, and 
it defines $\{ c_1 , \ldots c_N \}$ inside $L[(c_1 , \ldots c_N)]$. $\square$

\begin{lemma}\label{sacks-prop} {\bf (Sacks property)} Let $N<\omega$, and let 
$g$ be ${\mathbb P}(N)$-generic over $L$.
For each $f \colon \omega \rightarrow \omega$, $f \in L[a]$, there is some $g \in
L$ with domain $\omega$ such that for each $n<\omega$, $f(n) \in g(n)$ and \footnote{In what follows, the only thing that
will matter is that the bound on ${\rm Card}(g(n))$ only depends on $n$ and not on the particular $g$.}
${\rm Card}(g(n)) \leq (n+1) \cdot 2^{n+1}$.
\end{lemma}

{\em Proof.} Let $\tau \in L^{{\mathbb P}(N)}$, $\tau^g=f$. Let $(A_n \colon n<\omega) \in L$ be such that for each $n$, $A_n$ is a maximal antichain of 
${\vec T} \in {\mathbb P}(N)$ such that $\exists m<\omega \, {\vec T} \Vdash \tau({\check n})
= {\check m}$.
We may pick some $\xi < \omega_1$ such that $\bigcup \{ A_n \colon n<\omega \} \subset
{\mathbb P}_\xi(N)$ and $(A_n \colon n<\omega) = f_\xi$.

By Lemma \ref{5}, there are pairwise different $(m_1,s_1)$, $\ldots$, $(m_N,s_N)$
such that $$(A^{m_1}_{s_1,\xi}, \ldots , A^{m_N}_{s_N,\xi}) \in g.$$ 

Let $$n > {\rm max}( \{ N-1 , m_1 , \ldots , m_N ,
{\rm lh}(s_1) , \ldots , {\rm lh}(s_N) \} ){\rm . }$$ If $t_1 \supset s_1$, $\ldots$, $t_N \supset t_N$ are such that 
${\rm lh}(t_1) = \ldots = {\rm lh}(t_N) = n+1$,
then $(T^{m_1}_{t_1,\xi}, \ldots , T^{m_N}_{t_N,\xi}) \in d_\xi(n,N) \subset A_n$,
so that also $$\exists m<\omega \, (T^{m_1}_{t_1,\xi}, \ldots , T^{m_N}_{t_N,\xi}) \Vdash \tau({\check n}) = {\check m}.$$
Therefore, if we let
\begin{eqnarray*}
g(n) = \{ m<\omega \colon \exists t_1 \supset s_1, \ldots \exists t_N \supset t_N \, ({\rm lh}(t_1) = \ldots = {\rm lh}(t_N) = n+1 \wedge \\ (T^{m_1}_{t_1,\xi}, \ldots , T^{m_N}_{t_N,\xi}) \Vdash \tau({\check n}) = {\check m} ) \}{\rm , }
\end{eqnarray*}
then $(A^{m_1}_{s_1,\xi}, \ldots , A^{m_N}_{s_N,\xi}) \Vdash \tau({\check n})
\in (g(n)) {\check {}}$ ,  hence $f(n) \in g(n)$, and ${\rm Card}(g(n)) =
N \cdot 2^{n+1} \leq (n+1) \cdot 2^{n+1}$ for all but finitely many $n$.
$\square$

\section{The variant of the Cohen-Helpern-L\'evy model.}

Let us force with ${\mathbb P}(\omega)$ over $L$, 
and let $g$ be
a generic filter. Let $c_n$, $n<\omega$, denote the Jensen reals which $g$ adds.
Let us write $A = \{ c_n \colon n<\omega \}$ for the set of those Jensen reals. The
model 
\begin{eqnarray}\label{the_model!}
H = H(L) = {\sf HOD}^{L[g]}_{A \cup \{ A \}}
\end{eqnarray}
of all sets which inside
$L[g]$ are hereditarily definable from parameters in ${\rm OR} \cup
A \cup \{ A \}$ is the variant of the 
Cohen--Halpern--L\'evy model (over $L$) which we shall work with.
For the case of Jensen's original forcing this model was first considered in
\cite{enayat}.

For any finite
$a \subset A$, we write $L[a]$ for the model constructed from the finitely many
reals in $a$. 

\begin{lemma}\label{definability_2}
Inside $H$, $A$ is a (lightface) $\Pi^1_2$ set.
\end{lemma}

{\em Proof.} Let $\varphi(-)$ be the $\Pi^1_2$ formula from the proof of
Lemma \ref{definability}. If $H \models \varphi(x)$, $x \in L[a]$, $a \in [A]^{<\omega}$, then $L[a] \models \varphi(x)$ by Shoenfield, so $x \in a
\subset A$. On the other hand, if $c \in A$, then $L[c] \models \varphi(c)$
and hence $H \models \varphi(c)$ again by Shoenfield. $\square$

\bigskip
Fixing some G\"odelization of formulae (or some enumeration of all the rud
functions, resp.) at the outset, each $L[a]$, $a \in [A]^{<\omega}$, comes with a unique canonical global well--ordering $<_a$ of $L[a]$ by which we mean the one which is induced by the {\em 
natural} order of the elements of $a$ and the fixed G\"odelization device in the usual fashion.
The assignment $a \mapsto <_a$, $a \in [A]^{<\omega}$, is hence in $H$.\footnote{More
precisely, the ternary relation consisting of all $(a,x,y)$ such that $x <_a y$ is definable over $H$.}
This is a crucial fact.

Let us fix a bijection 
\begin{eqnarray}\label{e}
e \colon \omega \rightarrow \omega \times \omega{\rm , }
\end{eqnarray}
and let us write $((n)_0,(n)_1)= e(n)$.

We shall also make use the following. Cf.\ \cite[Lemma 1.2]{BSchWY}.

\begin{lemma}\label{no_wo}
(1) Let $a \in [A]^{<\omega}$ and $X \subset L[a]$, $X \in H$, say
$X \in {\sf HOD}^{L[g]}_{b \cup \{ A \}}$, where $b \supseteq a$, $b \in
[A]^{<\omega}$. Then $X \in L[b]$.

(2) There is no well--ordering of the reals in $H$.

(3) $A$ has no countable subset in $H$.

(4) $[A]^{<\omega}$ has no countable subset in $H$.
\end{lemma}

{\em Proof sketch.} (1) Every permutation $\pi \colon \omega \rightarrow \omega$ induces
an automorphism $e_\pi$ of ${\mathbb P}(\omega)$ by sending $p$ to
$q$, where $q(\pi(n))=p(n)$ for all $n<\omega$. It is clear that no $e_\pi$ moves
the canonical name for $A$, call it ${\dot A}$. Let us also write ${\dot c}_n$
for the canonical name for $c_n$, $n<\omega$. Now if $a$, and $b$ are as
in the statement of (1), say $b = \{ c_{n_1}, \ldots , c_{n_k} \}$,
if $p$, $q \in {\mathbb P}(\omega)$, if $\pi \upharpoonright \{ n_1 , \ldots ,
n_k \} = {\rm id}$, $p \upharpoonright \{ n_1 , \ldots ,
n_k \}$ is compatible with $q \upharpoonright \{ n_1 , \ldots ,
n_k \}$, and ${\rm supp}(\pi(p)) \cap {\rm supp}(q) \subseteq \{ n_1 , \ldots ,
n_k \}$, if $x \in L$, if $\alpha_1$, $\ldots$, $\alpha_m$
are ordinals, and if $\varphi$ is a formula, then 
\begin{eqnarray*}
p \Vdash_L^{{\mathbb P}(\omega)} \, \varphi({\check x},{\check \alpha}_1, 
\ldots {\check \alpha}_m, {\dot c}_{n_1} , \ldots {\dot c}_{n_k}, {\dot A}) & \Longleftrightarrow \\
\pi(p) \Vdash_L^{{\mathbb P}(\omega)} \, \varphi({\check x},{\check \alpha}_1, 
\ldots {\check \alpha}_m, {\dot c}_{n_1} , \ldots {\dot c}_{n_k}, {\dot A}) & {}
\end{eqnarray*}
and $\pi(p)$ is compatible with $q$,
so that the statement $\varphi({\check x},{\check \alpha}_1, 
\ldots {\check \alpha}_m, {\dot c}_{n_1} , \ldots {\dot c}_{n_k}, {\dot A})$
will be decided by conditions $p \in {\mathbb P}(\omega)$ with ${\rm supp}(p)
\subseteq  \{ n_1 , \ldots ,
n_k \}$. But every set in $L[b]$ is coded by a set of ordinals, so if $X$ is as in (1), this shows that $X \in L[b]$. 

(2) Every real is a subset of $L$. Hence by (1), if $L[g]$ had a well--ordering of the 
reals in ${\sf HOD}^{L[g]}_{a \cup \{ A \}}$, some $a \in [A]^{<\omega}$, then
every real of $H$ would be in $L[a]$, which is nonsense.  

(3) Assume that $f \colon \omega \rightarrow A$ is injective, $f \in H$.
Let $x \in {}^\omega \omega$ be defined by $x(n)=f((n)_0)((n)_1)$, so 
that $x \in H$. By (1), $x \in L[a]$ for some 
$a \in [A]^{<\omega}$. But then ${\rm ran}(f) \subset L[a]$,
which is nonsense, as there is some $n<\omega$ such that
$c_n \in {\rm ran}(f) \setminus a$. 

(4) This readily follows from (3). \hfill $\square$ (Lemma \ref{no_wo})

\bigskip
Let us recall another standard fact. 
\begin{eqnarray}\label{intersection}
\mbox{If $a$, $b \in [A]^{<\omega}$, then } L[a] \cap L[b] = L[a \cap b].
\end{eqnarray}
To see this,
let us assume without loss of generality that $a \setminus b \not= \emptyset \not=
b \setminus a$, and say $a \setminus b = \{ c_n \colon n \in I \}$ and
$b \setminus a = \{ c_n \colon n \in J \}$, where $I$ and $J$ are non--empty disjoint
finite subsets of $\omega$. Then 
$a \setminus b$ and $b \setminus a$ are mutually 
${\mathbb P}(I)$- and ${\mathbb P}(J)$-generic over $L[a \cap b]$. But then $L[a] \cap L[b] =
L[a \cap b][a \setminus b] \cap L[a \cap b][b \setminus a] = L[a \cap b]$,
cf.\ \cite[Problem 6.12]{book}.  

For any $a \in [A]^{<\omega}$, 
we write ${\mathbb R}_a = {\mathbb R} \cap L[a]$
and ${\mathbb R}^+_a = {\mathbb R}_a \setminus \bigcup \{ {\mathbb R}_b
\colon b \subsetneq a \}$. 
$({\mathbb R}^+_a
\colon a \in [A]^{<\omega})$ is a partition of ${\mathbb R}$: By Lemma
\ref{no_wo} (1),
\begin{eqnarray}\label{reals1}
{\mathbb R} \cap H = \bigcup \{ {\mathbb R}_a^+ \colon a \in 
[A]^{<\omega} \}{\rm , }
\end{eqnarray}
and ${\mathbb R}_a \cap {\mathbb R}_b = {\mathbb R}_{a \cap b}$
by (\ref{intersection}), so that
\begin{eqnarray}\label{reals2}
{\mathbb R}_a^+ \cap {\mathbb R}_b^+ = \emptyset \mbox{ for }
a, b \in [A]^{<\omega}, a \not= b.
\end{eqnarray}

For $x \in {\mathbb R}$, we shall also write $a(x)$ for the unique $a \in [A]^{<\omega}$
such that $x \in {\mathbb R}^+_a$, and we shall write $\#(x)={\rm Card}(a(x))$.

Adrian Mathias showed that in the original Cohen--Halpern--L\'evy model 
there is an definable function which assigns to each
$x$ an ordering $<_x$ such that $<_x$ is a well--ordering iff $x$ can be well--ordered,
cf.\ \cite[p.\ 182]{mathias}. The following is a special simple case of this,
adapted to the current model $H$.

\begin{lemma}\label{ctble_union} (A.\ Mathias)
In $H$, the union of countably many countable sets of reals is countable.
\end{lemma}

{\em Proof.} Let us work inside $H$. Let $(A_n \colon n<\omega)$ be such that
for each $n<\omega$, $A_n \subset {\mathbb R}$ and 
there exists some surjection $f \colon \omega \rightarrow 
A_n$. For each such pair $n$, $f$ let $y_{n,f} \in {}^\omega \omega$ be such
that $y_{n,f}(m)=f((m)_0)((m)_1)$. If $a \in [A]^{<\omega}$ and
$y_{n,f} \in {\mathbb R}_a$,
then $A_n \in L[a]$. By (\ref{intersection}), for each $n$ there is a unique $a_n \in [A]^{<\omega}$
such that $A_n \in L[a_n]$ and $b \supset a_n$ for each $b \in [A]^{<\omega}$
such that $A_n \in L[b]$. Notice that $A_n$ is also countable in $L[a_n]$.

Using the function $n \mapsto a_n$, an easy recursion yields a surjection $g \colon
\omega \rightarrow \bigcup \{ a_n \colon n<\omega \}$: first
enumerate the finitely many elements of $a_0$ according to their natural order,
then enumerate the finitely many elements of $a_1$ according to their natural order, etc.
As $A$ has no countable subset, $\bigcup \{ a_n \colon n<\omega \}$ must
be finite, say $a = \bigcup \{ a_n \colon n<\omega \} \in [A]^{<\omega}$.
But then $\{ A_n \colon n<\omega \} \subset L[a]$. (We don't claim
$(A_n \colon n<\omega) \in L[a]$.)

For each $n<\omega$, we may now let $f_n$ the $<_a$--least surjection $f \colon
\omega \rightarrow A_n$. Then $f(n) = f_{(n)_0}((n)_1)$ for $n<\omega$
defines a surjection from $\omega$ onto $\bigcup \{ A_n \colon n<\omega \}$, as desired.
\hfill $\square$ (Lemma \ref{ctble_union})

The following is not true in the original Cohen--Halpern--L\'evy model. Its proof
exploits the Sacks property, Lemma \ref{sacks-prop}. 

\begin{lemma}\label{null_set} (1) Let $M \in H$ be a null set in $H$.
There is then a $G_\delta$ null set $M'$ with $M' \supset M$ whose code is in $L$.

(2) Let $M \in H$ be a meager set in $H$.
There is then an $F_\sigma$ meager set $M'$ with $M' \supset M$ whose code is in $L$.
\end{lemma}

{\em Proof.} (1)
Let $M \in H$ be a null set in $H$. 

Let us work in $H$. Let $(\epsilon_n \colon n<\omega)$ be any sequence of positive reals.
Let $\bigcup_{s \in X} \, U_s \supset H$, where $X \subset {}^{<\omega} 2$ and
$\mu(\bigcup \{ U_s \colon s \in X \}) \leq \epsilon_0$.\footnote{Here, $\mu$ denotes Lebesge measure.}
Let $e \colon\omega \rightarrow X$ be onto. Let $(k_n \colon n<\omega)$ be defined by:
$k_n =$ the smallest $k$ (strictly bigger than $k_{n-1}$ if $n>0$) such that $\mu(\bigcup \{ U_s \colon s \in e
\mbox{''} \omega \setminus k \}) \leq \epsilon_{n}$.
Write $k_{-1}=0$. We then have that $\mu(\bigcup \{ U_s \colon s \in e \mbox{''} [k_{n-1},k_n) \}) \leq \epsilon_n$ for every $n<\omega$.

Now fix $\epsilon >0$. Let $$\epsilon_n = \frac{\epsilon}{n \cdot 2^{2n+2}}{\rm , }$$ and let $(k_n
\colon n<\omega)$ and $e \colon\omega \rightarrow {}^{< \omega} 2$ be such that
$\bigcup_{s \in X} \, U_s \supset H$ and $\mu(\bigcup \{ U_s \colon s \in e \mbox{''} [k_{n-1},k_n) \}) \leq \epsilon_n$ for every $n<\omega$.
We may now apply Lemma \ref{sacks-prop} inside $L[a]$ for some 
$a \in [A]^{<\omega}$ such that $\{ e , (k_n \colon n<\omega ) \} \subset L[a]$
and 
find a function $g \in L$ with domain $\omega$
such that for each $n<\omega$, $g(n)$ is a finite union $U_n$ of basic open sets
such that $\{ U_s \colon s \in e \mbox{''} [k_{n-1},k_n) \} \subset U_n$
and $\mu(U_n) \leq \frac{1}{2^{n+1}}$. But then ${\cal O} =
\bigcup \{ O_n \colon n<\omega \} \supset M$ is open, ${\cal O}$ is coded in $L$
(i.e., there is $Y \in L$, $Y \subset {}^{<\omega} 2$, with ${\cal O} =
\bigcup \{ U_s \colon s \in Y \}$), and $\mu({\cal O}) \leq \epsilon$.

We may hence for every $n<\omega$ let ${\cal O}_n$ be an open set with
${\cal O}_n \supset M$, $\mu({\cal O}_n) \leq \frac{1}{n+1}$, and whose
code in $L$ is  
$<_L$-least among all the codes giving such a set. Then $\bigcap \{ {\cal O}_n \colon n<\omega \}$
is a $G_\delta$ null set with code in $L$ and which covers $M$.

(2) Let $M \in H$ be a meager set in $H$, say $M = \bigcup \{ N_n \colon n<\omega \}$, where each $N_n$ is nowhere dense. 

Let us again work in $H$. It is easy to verify
that a set $P \subset {}^\omega 2$ is nowhere dense iff there is some $z \in {}^\omega 2$ and some strictly increasing $(k_n \colon n<\omega)$ such that for all $n<\omega$,
\begin{eqnarray}
\{ x \in {}^\omega 2 \colon x \upharpoonright [k_n,k_{n+1}) = z \upharpoonright [k_n,k_{n+1}) \} \cap P = \emptyset.
\end{eqnarray}
Look at $f \colon \omega \rightarrow \omega$, where $f(m)=k_{n+1}$ for the least $n$ with $m \leq k_n$. We may first 
apply Lemma \ref{sacks-prop} inside $L[a]$ for some 
$a \in [A]^{<\omega}$ such that $f \in L[a]$
and get a function $g \colon \omega \rightarrow \omega$, $g \in L$, such that $g(m) \geq f(m)$ for all $m<\omega$. Write $\ell_0=0$ and $\ell_{n+1} = g(\ell_n)$, so that
for each $n$ there is some $n'$ with 
\begin{eqnarray}
\ell_n \leq k_{n'} < k_{n'+1} \leq \ell_{n+1}.
\end{eqnarray}
Define $e \colon \omega \rightarrow \omega$ by $e(n) = \sum_{q=0}^n \, (q+1) \cdot
2^{q+1}$. We may now apply Lemma \ref{sacks-prop} inside $L[a]$ for some 
$a \in [A]^{<\omega}$ such that $f \in L[a]$
and get some $n \mapsto (z^n_i \colon i \leq (n+1) \cdot 2^{n+1})$ inside $L$ such that
for all $n$, $i$, $z^n_i \colon e(n) \rightarrow 2$, and for all $n$ there is some
$i$ with $z \upharpoonright e(n) = z^n_i$. From this we get some $z' \colon \omega
\rightarrow \omega$, $z' \in L$, such that for all $n$ 
there is some $n'$ with $z' \upharpoonright [\ell_{n'}, \ell_{n'+1}) = z
\upharpoonright [\ell_{n'}, \ell_{n'+1})$. But then, writing 
\begin{eqnarray}
D = \{ x \in {}^\omega 2 \colon \exists n \, x \upharpoonright [\ell_{e(n)},\ell_{e(n+1)}) = z' \upharpoonright [\ell_{e(n)},\ell_{e(n+1)}) \} {\rm , }
\end{eqnarray}
$D \in L$, and $D$ is open and dense.

We may hence for every $n<\omega$ let ${\cal O}_n$ be an open dense set with
${\cal O}_n \cap N_n = \emptyset$, whose
code in $L$ is  
$<_L$-least among all the codes giving such a set. Then $\bigcup \{ {}^\omega 2 \setminus {\cal O}_n \colon n<\omega \}$
is an $F_\sigma$ meager set with code in $L$ and which covers $M$.
$\square$ 

\begin{cor}\label{sierpinski_set}
In $H$, there is a $\Delta^1_2$ Sierpi\'nski set as well as a $\Delta^1_2$ Luzin set.
\end{cor}

{\em Proof.} There is a $\Delta^1_2$ Luzin set in $L$.
By Lemma \ref{null_set} (2), any such set is still a Luzin set in $H$.
The same is true with ``Luzin'' replaced by ``Sierpi\'nski'' and Lemma \ref{null_set} (2) replaced by Lemma \ref{null_set} (1).
$\square$

\begin{lemma}\label{bernstein_set}
In $H$, there is a $\Delta^1_3$ Bernstein set.
\end{lemma}

{\em Proof.} In this proof, let us think of reals as elements of the Cantor space
${}^\omega 2$. Let us work in $H$.

We let
\begin{eqnarray*}
B = \{ x \in {\mathbb R} \colon \exists \mbox{ even } n \, (2^n < \#(x) \leq 2^{n+1}) \} & \mbox{ and} \\
B' = \{ x \in {\mathbb R} \colon \exists \mbox{ odd } n \, (2^n < \#(x) \leq 2^{n+1}) \}. & {}
\end{eqnarray*}
Obviously, $B \cap B' = \emptyset$.

Let $P \subset {\mathbb R}$ be perfect. We aim to see that $P \cap B \not=
\emptyset \not= P \cap B'$.

Say $P = [T] = \{ x \in {}^\omega 2 \colon \forall n \, x \upharpoonright n \in T \}$,
where $T \subseteq {}^{<\omega} 2$ is a perfect tree. Modulo some fixed
natural bijection ${}^{<\omega} 2 \leftrightarrow \omega$, we may identify $T$ with a real.
By (\ref{reals1}),
we may pick some $a \in [A]^{<\omega}$ such that $T \in L[a]$.
Say ${\rm Card}(a) < 2^n$, where $n$ is even.

Let $b \in [A]^{2^{n+1}}$, $b \supset a$, and let $x \in {\mathbb R}^+_b$.
In particular, $\#(x)=2^{n+1}$.
It is easy to work in $L[b]$ and construct some $z \in [T]$ such that $x \leq_T z 
\oplus T$,\footnote{Here, $(x \oplus y)(2n)=x(n)$ and $(x \oplus y)(2n+1)=y(n)$, $n
<\omega$.} e.g., arrange that if $z \upharpoonright m$ is the $k^{\rm th}$ splitting
node of $T$ along $z$, where $k \leq m < \omega$, then $z(m)=0$ if $x(k)=0$
and $z(m)=1$ if $x(k)=1$.

If we had $\#(z) \leq 2^n$, then $\#(z \oplus T) \leq \#(z) + \#(T) < 2^n + 2^n =
2^{n+1}$, so that $\#(x) < 2^{n+1}$ by $x \leq_T z 
\oplus T$. Contradiction! Hence $\#(z) > 2^n$. By $z \in L[b]$, $\#(z) \leq 2^{n+1}$.
Therefore, $z \in P \cap B$. 

The same argument shows that $P \cap B' \not= \emptyset$. $B$ (and also $B'$) is
thus a Bernstein set. 

We have that $x \in B$ iff 
\begin{gather*}
\exists a \in [A]^{<\omega} \, \exists \mbox{ even } n \, \exists J_\alpha[a] \, \\ (x \in J_\alpha[a] \wedge 2^n < {\rm Card}(a) \leq 2^{n+1} \wedge 
\forall b \subsetneq a \, \forall J_\beta[b] x \notin J_\beta[b] ){\rm , }
\end{gather*}
which is true iff
\begin{gather*}
\forall a \in [A]^{<\omega} \, \forall J_\alpha[a] \, (x \in J_\alpha[a] \rightarrow 
\exists a' \subset a \, \exists \mbox{ even } n \, \exists J_{\alpha'}[a'] \, \\ (x \in J_{\alpha'}[a'] \wedge 2^n < {\rm Card}(a) \leq 2^{n+1} \wedge 
\forall b \subsetneq a' \, \forall J_\beta[b] x \notin J_\beta[b] )){\rm . }
\end{gather*}
By Lemma \ref{definability_2}, this shows that $B$ is $\Delta^1_3$.
$\square$ 

\bigskip
Recall that for any $a \in [A]^{<\omega}$, 
we write ${\mathbb R}_a = {\mathbb R} \cap L[a]$.
Let us now also write ${\mathbb R}_{<a}
=  {\rm span}( \bigcup 
\{ {\mathbb R}_{b} \colon b \subsetneq a \} )$, and 
${\mathbb R}^*_a = {\mathbb R}_a \setminus {\mathbb R}_{<a}$.
In particular, ${\mathbb R}_{<\emptyset} = \{ 0 \}$ by
our above convention that ${\rm span}(\emptyset)=\{ 0 \}$, and ${\mathbb R}^*_{\emptyset} = ({\mathbb R} \cap L)
\setminus \{ 0 \}$.

The proof of Claim \ref{claim1} below will show that
\begin{eqnarray}\label{reals3} 
{\mathbb R} \cap H = {\rm span}(\bigcup \{ {\mathbb R}^*_a \colon a \in [A]^{<\omega} \}).
\end{eqnarray}
Also, we have that ${\mathbb R}_a^* \subset {\mathbb R}_a^+$, so that by
(\ref{reals2}),
\begin{eqnarray}\label{reals4} {\mathbb R}_a^* \cap {\mathbb R}_b^* = \emptyset \mbox{ for }
a, b \in [A]^{<\omega}, a \not= b.
\end{eqnarray}

%

\begin{lemma}\label{main_thm}
In $H$, there is a $\Delta^1_3$ Hamel basis.
\end{lemma}

{\em Proof.} 
We call $X \subset {\mathbb R}^*_a$ {\em linearly independent over} ${\mathbb 
R}_{<a}$ iff whenever $$\sum_{n=1}^m \, q_n \cdot x_n \in 
{\mathbb R}_{<a}{\rm , }$$ where $m \in {\mathbb N}$, $m \geq 1$, and $q_n \in
{\mathbb Q}$ and $x_n \in X$ for all $n$, $1 \leq n \leq m$,
then $q_1 = \ldots = q_m = 0$. 
In other words, $X \subset {\mathbb R}^*_a$ is linearly independent over ${\mathbb 
R}_{<a}$ iff $${\rm span}(X) \cap {\mathbb R}_{<a} = \{ 0 \}.$$
We call $X \subset {\mathbb R}^*_a$ {\em maximal linearly independent over} ${\mathbb 
R}_{<a}$ iff $X$ is linearly independent over ${\mathbb 
R}_{<a}$ and no $Y \supsetneq X$, $Y \subset {\mathbb R}_a^*$ is still
linearly independent over ${\mathbb R}_{<a}$.
In particular, $X \subset {\mathbb R}^*_\emptyset = ({\mathbb R} \cap L) \setminus 
\{ 0 \}$ is linearly independent over ${\mathbb R}_{<\emptyset} = \{ 0 \}$ iff
$X$ is a Hamel basis for ${\mathbb R} \cap L$.
 
For any $a \in [A]^{<\omega}$, 
we let $b_a = \{ x_i^a \colon i<\theta^a \}$, some $\theta^a \leq \omega_1$, 
be the unique set such that
\begin{enumerate}
\item[(i)] for each $i < \theta^a$, $x_i^a$ is the $<_a$-least $x \in 
{\mathbb R}_a^*$ such that
$\{ x_j^a \colon j<i \} \cup \{ x \}$ is linearly independent over ${\mathbb R}_{<a}$,
and
\item[(ii)] $b_a$ is maximal linearly independent over ${\mathbb 
R}_{<a}$.
\end{enumerate}
By the above crucial fact, the function $a \mapsto b_a$ is well--defined 
and {\em exists inside} $H$. In particular, $$B = \bigcup \{ b_a \colon a \in [A]^{<\omega} \}$$
is an element of $H$.

We claim that $B$ is a Hamel basis for the reals of $H$, which will be established by Claims
\ref{claim1} and \ref{claim2}.

\begin{claim}\label{claim1} ${\mathbb R} \cap H \subset {\rm span}(B)$.
\end{claim}

{\em Proof} of Claim \ref{claim1}. Assume not, and let $n<\omega$ be the least
size of some $a \in [A]^{<\omega}$ such that ${\mathbb R}^*_a \setminus
{\rm span}(B) \not= \emptyset$. Pick $x \in {\mathbb R}^*_a \setminus
{\rm span}(B) \not= \emptyset$, where ${\rm Card}(a)=n$.

We must have $n>0$, as $b_\emptyset$ is a Hamel basis for the reals of $L$.
Then, by the maximality of $b_a$, while $b_a$ is 
linearly independent over ${\mathbb R}_{<a}$, $b_a \cup \{ x \}$ cannot be linearly independent
over ${\mathbb R}_{<a}$. This means that there are $q \in {\mathbb Q}$, $q \not= 0$, $m \in {\mathbb N}$, $m \geq 1$, and $q_n \in
{\mathbb Q} \setminus \{ 0 \}$ and $x_n \in b_a$ for all $n$, $1 \leq n \leq m$, such that
$$z = q \cdot x + \sum_{n=1}^m \, q_n \cdot x_n \in 
{\mathbb R}_{<a}.$$ By the definition of ${\mathbb R}_{<a}$ and the 
minimality of $n$, $z \in {\rm span}(\bigcup \{ b_c \colon c \subsetneq a \})$,
which then clearly implies that $x \in {\rm span}(\bigcup \{ b_c \colon c \subseteq a 
\}) \subset {\rm span}(B)$. 

This is a contradiction! \hfill $\square$ (Claim \ref{claim1})

\begin{claim}\label{claim2} $B$ is linearly independent.
\end{claim}

{\em Proof} of Claim \ref{claim2}. Assume not. This means that there are $1 \leq k < \omega$,
$a_i \in [A]^{<\omega}$ pairwise different,
$m_i \in {\mathbb N}$, $m_i \geq 1$ for $1 \leq i \leq k$, and $q_n^i \in
{\mathbb Q} \setminus \{ 0 \}$ and $x_n^i \in b_{a_i}$ for all 
$i$ and $n$ with $1 \leq i \leq k$ and $1 \leq n \leq m_i$ such that
\begin{eqnarray}\label{nobasis}
\sum_{n=1}^{m_1} \, q_n^1 \cdot x_n^1 + \ldots + 
\sum_{n=1}^{m_k} \, q_n^k \cdot x_n^k = 0.
\end{eqnarray} 
By the properties of $b_{a_i}$,
$\sum_{n=1}^{m_i} \, q_n^i \cdot x_n^i \in {\mathbb R}^*_{a_i}$,
so that (\ref{nobasis}) buys us that there are $z_i \in {\mathbb R}_{a_i}^*$,
$z_i \not= 0$,
$1 \leq i \leq k$,
such that 
\begin{eqnarray}\label{nobasis2}
z_1 + \ldots + z_k = 0.
\end{eqnarray}
There must be some $i$ such that there is no $j$ with $a_j \supsetneq a_i$,
which implies that $a_j \cap a_i \subsetneq a_i$ for all $j \not= i$. Let us assume without
loss of generality that $a_j \cap a_1 \subsetneq a_1$ for all $j$, $1<j \leq k$.

Let $a_1 = \{ c_\ell \colon \ell \in I \}$, where $I \in [\omega]^{<\omega}$,
and let $a_j \cap a_1 = \{ c_\ell \colon \ell \in I_j \}$, where $I_j \subsetneq I$,
for $1<j \leq l$. 

In what follows, a {\em nice name} $\tau$ for a real is a name of the
form 
\begin{eqnarray}\label{nicename}
\tau = \bigcup_{n,m<\omega} \, \{ (n,m)^\vee \} \times A_{n,m}{\rm , }
\end{eqnarray}
where each $A_{n,m}$ is a maximal antichain of conditions of the forcing in
question deciding that
$\tau({\check n})={\check m}$.

We have that $z_1$ is ${\mathbb P}(I)$--generic over $L$,
so that we may pick a nice name $\tau_1 \in L^{{\mathbb P}(I)}$
for $z_1$ with $(\tau_1)^{g \upharpoonright I}=z_1$.
Similarly, for $1 < j \leq k$,
$z_j$ is ${\mathbb P}(I_j)$--generic over
$L[g \upharpoonright (\omega \setminus I)]$, so that 
we may pick a nice name $\tau_j \in L[g \upharpoonright (\omega \setminus I)]^{{\mathbb P}(I_j)}$
for $z_j$ with $(\tau_j)^{g \upharpoonright I_j}=z_j$. We may construe each $\tau_j$,
$1<j \leq k$, as a name in $L[g \upharpoonright (\omega \setminus I)]^{{\mathbb P}(I)}$ by replacing each $p \colon I_j \rightarrow {\mathbb P}$ in an antichain
as in (\ref{nicename}) by $p' \colon I \rightarrow {\mathbb P}$, where
$p'(\ell)=p(\ell)$ for $\ell \in I_j$ and $p'(\ell) = \emptyset$ otherwise.
Let $p \in g \upharpoonright I$ be such that $$p \Vdash_{L[g \upharpoonright (\omega \setminus I)]}^{{\mathbb P}(I)} \tau_1+ \tau_2 + \ldots + \tau_k = 0.$$

We now have that inside $L[g \upharpoonright (\omega \setminus I)]$, there are
nice ${{\mathbb P}(I)}$--names $\tau_j'$, $1<j \leq k$ (namey, $\tau_j$, $1<j \leq k$),
such that still inside $L[g \upharpoonright (\omega \setminus I)]$  
\begin{enumerate}
\item[(1)] $p \Vdash^{{\mathbb P}(I)} \tau_1+ \tau_2' + \ldots + \tau_k' = 0$, and
\item[(2)] for all $j$, $1<j \leq k$ and for all $p$ in one of the antichains
of the nice name $\tau_j'$, ${\rm supp}(p) \subseteq I_j$.
\end{enumerate}
By Lemma \ref{4}, the nice names $\tau_1$, $\tau_2'$, $\ldots$, $\tau_k'$ may be
coded by reals, and both (1) and (2) are arithmetic in such real codes for $\tau_1$,$\tau_2'$, $\ldots$, $\tau_k'$,
so that by $\tau_1 \in L^{{\mathbb P}(I)}$ and 
$\Sigma^1_1$--absoluteness between $L$ and $L[g \upharpoonright (\omega \setminus I)]$ there are 
inside $L$ nice ${{\mathbb P}(I)}$--names $\tau_j'$, $1<j \leq k$, such that in $L$,
(1) and (2) hold true. But then, writing $z_j' = (\tau_j')^{g \upharpoonright I}$,
we have by (2) that $z_j' \in {\mathbb R}_{I_j}$ for $1<j \leq k$, and $z_1 +
z_2' + \ldots + z_k' =0$ by (1). But then $z_1 \in {\mathbb R}^*_I \cap 
{\mathbb R}_{<I}$, which is absurd. \hfill $\square$ (Claim \ref{claim2})

We now have that $x \in B$ iff
\begin{gather*}
\exists a \in [A]^{<\omega} \, \exists J_\alpha[a] \, \exists (x_i \colon i \leq \theta) \in J_\alpha[a] \, \exists X \subset \theta+1 \, (\mbox{ the } x_i \mbox{ enumerate the first } \\ \theta+1 \mbox{ reals in } J_\alpha[a] \mbox{ acc.\ to } <_a \wedge \, \theta \in X \wedge x = x_\theta \wedge \\
\forall i \in \theta \setminus X \exists J_\beta[a] \, J_\beta[a] \models 
\{ x_j \colon j \in X \cap i \} \cup \{ x_i \} \mbox{ is not linearly independent over } {\mathbb R}_{<a} \wedge \\
\forall i \in X \, \forall J_\beta[a] \, J_\beta[a] \models 
\{ x_j \colon j \in X \cap i \} \cup \{ x_i \} \mbox{ is linearly independent over } {\mathbb R}_{<a}))
{\rm , }
\end{gather*}
which is true iff
\begin{gather*}
\forall a \in [A]^{<\omega} \, \forall J_\alpha[a] \, \forall (x_i \colon i \leq \theta) \in J_\alpha[a] \, \forall X \subset \theta+1 \, (( \mbox{ the } x_i \mbox{ enumerate the first } \\ \theta+1 \mbox{ reals in } J_\alpha[a] \mbox{ acc.\ to } <_a \wedge \, x = x_\theta \wedge \\
\forall i \in (\theta+1) \setminus X \exists J_\beta[a] \, J_\beta[a] \models 
\{ x_j \colon j \in X \cap i \} \cup \{ x_i \} \mbox{ is not linearly independent over } {\mathbb R}_{<a} \wedge \\
\forall i \in X \, \forall J_\beta[a] \, J_\beta[a] \models 
\{ x_j \colon j \in X \cap i \} \cup \{ x_i \} \mbox{ is linearly independent over } {\mathbb R}_{<a}) \rightarrow \\
\theta \in X)
{\rm . }
\end{gather*}

By Lemma \ref{definability_2}, this shows that $B$ is $\Delta^1_3$. $\square$ 

\section{Open questions.}

We finish by stating some open problems.

\medskip
(1) Is there a model of ${\sf ZF}$ plus $\lnot {\sf AC}_\omega({\mathbb R})$
where there are sets as in (a)-(d) of Theorem \ref{main-thm} of lower projective complexity?

\medskip
(2) Does the model $H$ from (\ref{the_model!}) on p.\ \pageref{the_model!}
have a Burstin basis? An affirmative answer along the lines of the argument from \cite{brendle_et_al} would require us to show that 
\begin{eqnarray}\label{s_null}
{\mathbb R}_{<a} \in (s^0)^{L[a]} \mbox{ for all } a \in [A]^{<\omega}
{\rm , }
\end{eqnarray}
where $s^0$ denotes the Marczewski ideal. We don't know if (\ref{s_null}) is true, though, we don't even know if 
\begin{eqnarray}\label{s_null_2}
{\mathbb R}_{a}^* \not= \emptyset \mbox{ for all } a \in [A]^{<\omega}. 
\end{eqnarray}
L.\ Wu and L.\ Yu have recently shown that (\ref{s_null_2}) is true for ${\rm Card}(a)=2$, but
it is not known if (\ref{s_null_2}) holds true for ${\rm Card}(a)=3$. The second author has shown that if $A$ is a countable set of Cohen reals over $L$ (or, for that matter, any countable
set of dominating reals over $L$), then (\ref{s_null_2}) is true for $a \in [A]^{<\omega}$ of arbitrary size, i.e., that (\ref{s_null_2}) holds true for ${\mathbb R}_{a}^*$ as being defied in \cite{BSchWY}.

\medskip
(3) Does the model $H$ from (\ref{the_model!}) on p.\ \pageref{the_model!}
have a Mazurkiewicz set?

\medskip
We may force with the forcings from \cite{brendle_et_al} and \cite{mazurkiewicz}
to add a Burstin basis and a Mazurkiewicz set, respectively, over $H$ (without adding any reals), but then those
sets won't be definable in that extension.

\end{document}